\documentclass{amsart}[10pt,a4paper]
\usepackage[latin1]{inputenc}
\usepackage[english]{babel}
\usepackage{geometry}
\usepackage{array}
\usepackage{enumerate}
\usepackage{amsmath, amsfonts, amscd, amssymb, amsthm, xypic}

\usepackage{todonotes}

\newcommand{\sM}{s\mathcal{M}}

\newcommand{\D}{\mathcal{D}}

\renewcommand{\P}{\mathbb{P}}

\renewcommand{\c}{\mathcal{C}}

\newcommand{\n}{\mathfrak{n}}

\newcommand{\Z}{\mathbb{Z}}

\renewcommand{\S}{\mathcal{S}}
\renewcommand{\O}{\mathcal{O}}

\DeclareMathOperator{\Aut}{Aut}
\DeclareMathOperator{\Ber}{Ber}

\DeclareMathOperator{\rk}{rk}
\DeclareMathOperator{\red}{bos}
\DeclareMathOperator{\spli}{split}

\newenvironment{sis}{\left\{\begin{aligned}}{\end{aligned}\right.}

\theoremstyle{plain}
\newtheorem{theorem}{Theorem}[section]

\newtheorem{definition}[theorem]{Definition}
\newtheorem{criterion}[theorem]{Criterion}

\theoremstyle{definition}
\newtheorem{lemma}[theorem]{Lemma}

\newtheorem{remark}{Remark}[section]

\title{Pluri-Canonical Models of Supersymmetric Curves}
\date{}
\author[G. Codogni]{Giulio Codogni}
\address{Dipartimento di Matematica e Fisica, Universit\`{a} Roma Tre, Largo San Leonardo Murialdo, 1
00146, Roma, Italy.}
\email{codogni@mat.uniroma3.it}

\begin{document}

\maketitle

\begin{abstract}

This paper is about pluri-canonical models of supersymmetric (susy) curves. Susy curves are generalisations of Riemann surfaces in the realm of super geometry. Their moduli space is a key object in supersymmetric string theory. We study the pluri-canonical models of a susy curve, and we make some considerations about Hilbert schemes and moduli spaces of susy curves.

\end{abstract}

\begin{section}*{Introduction}
 
 This paper is about supersymmetric (susy) curves and their pluri-canonical models. Susy curves are generalisations of Riemann surfaces in the realm of super geometry; because of this, they are also called super Riemann surfaces. Their moduli space is a key object in supersymmetric string theory. 

\smallskip

Our goal is to construct the pluri-canonical model of a supersymmetric curve. In super geometry, the generalisation of the canonical bundle is the Berezinian bundle. Our result shows that an appropriate power of the Berezinian line bundle embeds all genus $g$ susy curves in a projective space of fixed dimension. The embedded susy curve is called the pluri-canonical model of the curve.

The prefix ``pluri" indicates that we are using powers of the Berezinian bundle.  As explained in Remark \ref{canonical}, a canonical model of a susy curve does not exist.

\begin{theorem}[=Theorem \ref{CanMod}, Pluri-Canonical Model]
Let $$f\colon \c\to B$$ be a genus $g\geq 2$ susy curve over a super-manifold $B$. Let $\Ber(\c)$ be the Berezinian of the relative cotangent bundle. Fix an integer $\nu\geq 3$. Then, $f_*(\Ber(\c)^{\otimes \nu})$ is a vector bundle on $B$.

Moreover, if the genus $g$ of $\c$ is at least $4$, the line bundle $\Ber(\c)^{\otimes \nu}$ is relatively very ample; in other words, it embeds $\c$ into $\P (f_*(\Ber(\c)^{\otimes \nu}))$. If $g=3$ or $g=2$, we need respectively $\nu\geq4$ or $\nu \geq 5$ for the very ampleness of $\Ber(\c)^{\otimes \nu}$. 
\end{theorem}
In Theorem  \ref{CanMod}, the reader can find also the rank of $f_*(\Ber(\c)^{\otimes \nu})$. The relevant definitions are given in Sections \ref{preliminaries} and \ref{definition}. The argument mimics the proof of the classical Kodaira's embedding theorem, see also \cite{ProjEmbed}. In Section \ref{preliminaries} we prove some general criteria to show that the push-forward of a line bundle is a vector bundle, and to check its very ampleness.

\smallskip

Our result is a prerequisites for the study of moduli space of susy curves using Hilbert schemes; we make some speculations in this direction in Section \ref{hilb}.

\smallskip

For the reader convenience, let us finish this introduction with a (non-comprehensive) review of the literature about the topics of this paper. General references about super geometry are \cite{Manin1}, \cite{Deligne}, \cite{Fioresi}, \cite{Witten1} and the first section of \cite{WD}; let us mention also the slightly more technical \cite{ProjEmbed}. Two important references about supersymmetric curves are \cite{Manin2} and \cite{Witten2}; other sources are \cite{Ber}, \cite{Rabin}, \cite{FK} and \cite{K}. Reference about moduli of susy curves are \cite{LeBrun} and \cite{spain}, two recent important contributions are \cite{WD} and \cite{WD2}. Recently, P. Deligne posted on his webpage a manuscript letter that he sent to Manin in 1987: this document contains a number of deep results probably still unexploited  (\cite{DeligneLetter}).

In the forthcoming paper \cite{Filippo}, we will present three different versions of the moduli space of susy curves, and  discuss the super-analogue of the period map.

\begin{subsection}*{Acknowledgements}
This paper is heavily influenced by the work of Edward Witten and some private conversations with him. We also had the pleasure of discussing these topics with T. Covolo, R. Donagi, S. Kwok, R. Fioresi and F. Viviani.

The author is funded by the FIRB 2012 ``Moduli spaces and their applications".  He acknowledges the support by the Simon Centre for his participation to the ``Supermoduli Workshop", and by the ``De Giorgi Center" for his participation to the intensive research period ``Perspectives in Lie Theory". 

\end{subsection}
\end{section}

\begin{section}{Preliminaries on line bundles and push-forward}\label{preliminaries}
First, let us recall the definition of super-manifold. Let $(X,\O_X)$ be a $\mathbb{Z}_2$-graded commutative ringed space, i.e. $X$ is a topological space and $\O_X$ is a sheaf of super commutative algebras. Let 
$\n_X\subset \O_X$ (or simply $\n$) be the ideal sheaf generated by odd elements. We say that $(X,\O_X)$ (or simply $X$) is a smooth super-variety\footnote{We refrain from defining more general super-varieties since in this paper we will only deal with smooth super-varieties.}, or \emph{super-manifold}, if the following conditions hold:
\begin{enumerate}
\item $X^{\red}:=(X,\O_X/\n_X)$ is a complex manifold;
\item the sheaf $\n_X/\n_X^2$ is a locally free sheaf of $\O_X/\n_X$-modules;
\item  $\O_X$ is locally  isomorphic to the $\Z/2\Z$-graded exterior algebra $\Lambda^{\bullet}(\n_X/\n_X^2)$ over $\n_X/\n_X^2$.
\end{enumerate}
The dimension of a connected super-manifold, written as $n|m$, is a pair of numbers: the first, the even dimension, is the dimension of $X^{\red}$; the second, the odd dimension, is the rank of $\n_X/\n_X^2$.  
Note that $n|0$ super-manifolds are ordinary manifolds. For any super-manifold $X$, there is a natural  closed embedding of $X^{\red}\hookrightarrow X$ which is an identity on the underlying topological spaces and it is the quotient map $\O_X\twoheadrightarrow \O_X/\n_X$ on the sheaf of functions. The manifold $X^{\red}$ is sometime called the \emph{bosonic reduction}, or just the reduction when no confusion arise, of $X$.

\medskip

It is possible to construct a super-manifold out of a complex manifold $Y$ and a vector bundle $\pi \colon E\to Y$: the topological space is $Y$ and the structure sheaf is $\O_Y\otimes \bigwedge^{\bullet} E$. We denote this new super-manifold by $Y_E$. In this case, $Y_E^{\red}=Y$ and $\n_{Y_E}/\n_{Y_E}^2=E$. The projection $\pi$ induces a section of the natural inclusion $Y \hookrightarrow Y_E$. 

A super-manifold $X$ is called \emph{split} if it is isomorphic to to $Y_E$, where $Y$ has to be isomorphic to $X^{\red}$ and $E$ to $\n_{X}/\n_{X}^2$. Any $n|1$ dimensional super manifold is split, see e.g. \cite[Corollary 2.2]{WD}; on the other hand, super-manifolds with odd dimension strictly bigger than one are not, in general, split. 

A \emph{vector bundle} $E$ on a super-manifold $X$ is a locally free sheaf of $\O_X$-modules. Let $U$ be an open subset over which $E$ is free. On $U$, $E$ is of the form $A^p|A^q$, where $A=\O_X(U)$. We define the rank of $E$ to be $p|q$. The definition does not depend on $U$. More details can be found in \cite[Section IV.9]{Manin1} or \cite[Appendix B3]{Fioresi}. We define $E^{\red}$ to be the pull-back of $E$ to $X^{\red}$. More concretely, $E^{\red}$ is $E/\n E$. The bundle $E^{\red}$ is a $\mathbb{Z}_2$-graded vector bundle on $X^{\red}$; its graded rank equals the rank of $E$. A \emph{line bundle} $E$ on $X$ is a vector bundle of rank either $1|0$ or $0|1$; in particular, $E^{\red}$ is a $\mathbb{Z}_2$ graded line bundle on $X^{\red}$. 
\medskip

In this paper, we are interested in the relative setting. Let $f\colon X\to B$ be a morphisms between smooth manifolds and let $f^{\red}:X^{\red}\to B^{\red}$ be the induced morphism on the corresponding reduced manifolds. We say that $f$ is proper if $f^{\red}$ is proper. We say that $f$ is smooth if its differential is surjective, i.e. the morphism is a submersion, see \cite[Section 5.2]{Fioresi}. The definition of flatness is as in the classical case, see \cite[Section 1.5]{Vaintrob} and \cite[Section 2.1.3]{SuperStacks}; in particular, a submersion is flat. Through all this paper, we assume that $f$ is proper and smooth and that the relative dimension is $1|1$. 

We define the split model $X^{\spli}$ of $X$ to be the cartesian product of $B^{\red}$ and $X$; let us draw the resulting commutative diagram.
\begin{displaymath}
\begin{array}{ccccc}
X^{\red} &   \xrightarrow{j} & X^{\spli}       & \xrightarrow{\iota} &     X\\
              &       \searrow^{f^{\red}}     & \downarrow^{f^{\spli}}  &                              & \downarrow^f \\
              &                          & B^{\red}        &  \xrightarrow{i}& B

\end{array}
\end{displaymath}

The idea is that $X^{\spli}$ is still a super-manifold in between $X$ and $X^{\red}$, where we have killed all the odd functions coming from the base $B$ - the odd moduli - but we still keep the odd functions living on the fibres of $f$. When the relative dimension is $1|1$, as in this paper, $X^{\spli}$ is really a split manifold; this false in a more general setting. Being split, we can write $X^{\spli}=C_L$, where $C=X^{\red}$ is a classical curve over $B^{\red}$ and $L$ is a line bundle on $C$. Let $\pi\colon X^{\spli}\to X^{\red}$ be a section of $j$. Such a $\pi$ exists because $X^{\spli}$ is a split variety; this splitting is unique up to a scalar, but we do not need and we do not prove this fact.

Being $X^{\spli}\cong C_L$, we have $\pi_*\O_{X^{\spli}}=\O_{X^{\red}}\oplus L$. For a line bundle $E$ on $X$, we have
$$
\pi_*\iota^*E=E^{\red}\oplus (E^{\red}\otimes L)
$$
In other words, this is a splitting of $E/\n_BE$ as $\O_{X^{\red}}$-module. Remark also that $\pi \circ f^{\red}=f^{\spli}$. The following criterion is a generalisation  of \cite[Section 2.6]{Ber}.

\begin{criterion}\label{loc_free}
Let $f\colon X \to B$ be a proper smooth morphism of relative dimension $1|1$ between super manifolds; let $E$ be a line bundle on $X$. Write $X^{\spli}=C_L$. Then, the sheaf $f_*E$ is a vector bundle on $B$ if the two following conditions hold
$$\begin{sis}
& f^{\red}_*E^{\red} \textrm{ and } f^{\red}_*(E^{\red}\otimes L)   \; \textrm{are vector bundles on} \; B^{\red};\\
& R^1f^{\red}_*E^{\red}=R^1f^{\red}_*(E^{\red}\otimes L)=0.\\
\end{sis}$$

Moreover, under these hypotheses, the rank of $f_*E$ is equal to  $\rk(f^{\red}_*E^{\red}) \, | \,  \rk(f^{\red}_* E^{\red}\otimes L)$ if $E$ is of rank $1|0$, and to $\rk(f^{\red}_*E^{\red}\otimes L) \, | \,  \rk(f^{\red}_* E^{\red})$ if $E$ is of rank $0|1$.
\end{criterion}
\begin{proof}
We work at a local ring of a point of $B$; recall that Nakayama's Lemma and all its consequences hold in the super setting, cf. \cite[Appendix B3]{Fioresi}. We denote by $\n_B$ we denote the ideal sheaf on $B$ generated by odd elements. We also fix a section $\pi \colon X^{\spli}\to X^{\red}$ of $j\colon X^{\red}\to X^{\spli}$.

We are going to show by induction on $l$ that, for every $l$, the sheaf $f_*(E\otimes_{\O_X}f^* \O_B/\n_B^l)$ is a locally free sheaf of $\O_B/\n_B^l$-modules of rank either $\rk(f^{\red}_*E^{\red}) \, | \,  \rk(f^{\red}_* E^{\red}\otimes L)$ if $E$ is of rank $1|0$, or $\rk(f^{\red}_*E^{\red}\otimes L )\, | \,  \rk(f^{\red}_* E^{\red})$ if $E$ is of rank $0|1$. This is enough to conclude because, for $l$ big enough, the ideal sheaf $\n_B^l$ is trivial.

\smallskip

We first take $l=1$. We have 
\begin{equation}\label{first_step}
f_*(E\otimes_{\O_X} f^*\O_B/\n_B)=f_*\iota_*\iota^*E=i_*f^{\spli}_*\iota^*E= i_*f^{\red}_*\pi_*\iota^*E=i_*(f^{\red}_*E^{\red} \oplus f^{\red}_*(E^{\red}\otimes L ))
\end{equation}

The right hand side is a locally free sheaf of $\O_{B^{\red}}$ modules of the requested rank by hypothesis.

\smallskip

We now assume the statement for $l$ and we prove it for $l+1$. Consider the exact sequence of sheaves on $X$
$$
0\to f^*(\n_B^l/\n_B^{l+1})\to f^*(\O_B/\n_B^{l+1})\to f^*(\O_B/\n_B^l)\to 0
$$
We tensor by $E$; since $E$ is flat on $X$, the sequence remains exact, so we obtain 
$$
0\to f^*(\n_B^l/\n_B^{l+1})\otimes_{\O_X}E\to f^*(\O_B/\n_B^{l+1})\otimes_{\O_X}E\to f^*(\O_B/\n_B^l)\otimes_{\O_X}E\to 0
$$

Remark that the action of $\n_B$ on $\n_B^l/\n_B^{l+1}$ is trivial, so we can write
$$ f^*(\n_B^l/\n_B^{l+1})\otimes_{\O_X}E= \iota_*  ((f^{\spli})^*(\n_B^l/\n_B^{l+1})\otimes_{\O_{X^{\spli}}}\iota^*E )$$

Here, we are viewing $\n_B^l/\n_B^{l+1}$ as a locally free sheaf of $\O_{B^{\red}}$-modules. We now apply $Rf_*$. Since both $\iota$ and $i$ are affine, we have 
$$(Rf_*)(\iota_*( \, - \,))=R(f\circ \iota)_*( \, - \,)=R(i\circ f^{\spli})_*( \, - \,)  = i_*R((f^{\spli})^*)(\, - \,)  $$

Since $\n_B^l/\n_B^{l+1}$ is locally free on $B^{\red}$, we can apply the projection formula to $Rf^{\spli}$. A part of the resulting long exact sequence is
\begin{align*}
\begin{split}
0\to i_*( \n_B^l/\n_B^{l+1}\otimes_{\O_{B^{\red}}}f^{\spli}_*\iota^*E ) \to  f_*(f^*\O_B/\n_B^{l+1} \otimes_{\O_X}  E )  \to f_*(f^*\O_B/\n_B^l \otimes_{\O_X} E) \to \\
\to i_*( \n_B^l/\n_B^{l+1}\otimes_{\O_{B^{\red}}}(R^1f^{\spli}_*)\iota^*E )\,.
\end{split}
\end{align*}

Note that the last term vanishes:

\begin{align*}
\begin{split}
(R^1f^{\spli}_*)\iota^*E&= (R^1f^{\red}_*\pi_*)\iota^*E= (R^1f^{\red}_*)\pi_*\iota^*E=\\
& R^1f^{\red}_*E^{\red}\oplus R^1f^{\red}_*(E^{\red}\otimes L)=0 \,,
\end{split}
\end{align*}
where we used that $\pi$ is an affine maps, and the vanishing assumed in the hypotheses of the criterion.

We are now left with the short exact sequence
\begin{equation}\label{eq1}
0\to i_*( \n_B^l/\n_B^{l+1}\otimes_{\O_{B^{\red}}}f^{\spli}_*\iota^*E ) \to  f_*(f^*\O_B/\n_B^{l+1} \otimes_{\O_X}  E )  \to f_*(f^*\O_B/\n_B^l \otimes_{\O_X} E) \to 0
\end{equation}
As discussed above, we have
$$
 f^{\spli}_*\iota^*E=f^{\red}_*E^{\red} \oplus f^{\red}_*(E^{\red}\otimes L) 
 $$
so, the first term of the sequence (\ref{eq1}) is a locally free sheaf of $\O_B/\n_B$-modules by hypothesis; the last term is, by induction, a locally free sheaf of $O_B/\n_B^l$-modules of the requested rank. One can now apply a Nakayama-type argument to show that the central term is a locally-free sheaf of $\O_B/\n_B^{l+1}$-modules of the requested rank; see for example the proof of \cite[Tag 051H]{stacks-project} with $R=\O_B/\n_B^{l+1}$ and $I=\n_B^l$. More explicitly, one takes a basis for the last term of the sequence (\ref{eq1}); lift these elements to the central term, and shows that they generate it because of Nakayama's Lemma. Then, arguing again as in the proof of \cite[Tag 051H]{stacks-project}, one shows that, since $\n_B$ is nilpotent, the central term is locally free.

 \end{proof}

From now on, we assume that $f_*L$ is locally free. The line bundle $L$ gives a morphism relative to $B$ from $X$ to $\P f_*L$; we say that $L$ is $f$-\emph{very ample} if this morphism is an embedding. Following the proof of  \cite[Theorem 1]{ProjEmbed}, we have the following criterion. 

\begin{criterion}\label{very_ample}
Let $f\colon X \to B$ be a proper smooth morphism of relative dimension $1|1$ between super manifolds; write $X^{\spli}=C_L$. Let $E$ be a line bundle of rank $1|0$ on $X$. Assume that $f_*E$ is a locally free sheaf. Then $L$ is $f$-very ample if, for every pair of points $x$ and $y$ in $X^{\red}$ such that $f^{\red}(x)=f^{\red}(y)$ we have
 $$\begin{sis}
& R^1f^{\red}_* (E^{\red}\otimes I_x\otimes I_y)=0, \\
& R^1f^{\red}_*  ( E^{\red}\otimes L \otimes I_x)=0,   \\
\end{sis}$$
where $I_x\subset \O_{X^{\red}}$ (resp. $I_y$) is the ideal sheaf of $x$ (resp.  $y$). If $E$ is of rank $0|1$, the same statement holds replacing in the second condition $E^{\red}\otimes L$ with $E^{\red}$.
\end{criterion}
\begin{proof}
We have to show that the map given by $E$ embeds $X$ into $\P(f_*E)$. The first condition, by the classical proof of the Kodaira embedding theorem, guarantees that $E$ embeds $X^{\red}$ into $\P(f^{\red}_*E^{\red})$. We now have to verify the statement in the odd directions. The second condition means that the odd differential is injective in the vertical direction, and this is enough to conclude. 

Let us be more explicit;  to fix the notation, assume that the rank of $E$ is $1|0$. We also shrink $B$ around $f(x)$, so that $f_*E$ is free. The second vanishing in the hypothesis guarantees that the morphism
$$
f_*^{\red}(E^{\red}\otimes L) \to f_*^{\red}(E^{\red}_x\otimes L_x)
$$
is surjective, where the sub-index $x$ means the stalk. This means that there exists a global section $\bar{s}$ of $E^{\red}\otimes L$ which does not vanish at $x$. 

Using equation (\ref{first_step}), we can identify $\bar{s}$ with an \emph{odd} global section of $\iota^* E=f^*\O_B/\n_B\otimes E$. Using sequence (\ref{eq1}) enough times, we can show that there exists an odd global section $s$ of $E$ which restricts to $\bar{s}$ modulo $\n_B$. In particular, $s$ is not zero at $x$.

If we trivialise $E$ around $x$ on $X$, the global section $s$ restrict to a vertical odd co-ordinate at $x$.  (By vertical odd co-ordinate, we mean an odd co-ordinate which is not zero modulo $\n_B$.) This is enough to conclude that $E$ gives an embedding at $x$. More details in the proof of \cite[Theorem 1]{ProjEmbed}.

\end{proof}

\end{section}
\begin{section}{Pluri-canonical models of susy curves}

\begin{subsection}{Definition of susy curves}\label{definition}
We give the definition of susy curve over a smooth base $B$; the reason for this set up is twofold: it is suitable for moduli theory; and a curve defined over a point is always split, so there is not much super-geometry going on.

\begin{definition}[Susy Curves]

A genus $g$ susy curve over a smooth base $B$ is a proper smooth morphism
\[f \colon \c \to B\] 
between super manifolds together with a susy structure $\D$ satisfying the following conditions. The morphism $f$ has relative dimension is $1|1$ and the fibres are, as topological spaces, genus $g$ Riemann surfaces. The susy structure $\D$ is a sub bundle
\[\D \hookrightarrow Ker(df) \hookrightarrow T\c\]
which is of rank $0|1$, and such that $\D$ and $[\D,\D]$ span $Ker(df)$. 
\end{definition}
Smooth morphisms are discussed in \cite[Section 5.2]{Fioresi}, where they are called submersions. Recall that the coordinates on the base $B$ are sometime called the moduli of the curve.

\end{subsection}

\begin{subsection}{Pluri-Canonical model}
Let 
$$ f\colon \c \to B $$
be a susy curve; write $\Ber(\c):=\Ber(T_f\c^{\vee})$ for the Berezinian of the relative co-tangent bundle: this is rank $0|1$ line bundle. Since $\D$ and $[\D,\D]\cong \D^2$ spans $T_f\c$, we have an exact sequence
$$ 0\to \D^{-2} \to T_f\c^{\vee}\to \D^{-1}\to 0   $$
The Berezinian of this sequence gives an isomorphism $\Ber(\c)=\D^{-1}$; cf. \cite[Equation 2.28]{Witten2}. 

\begin{lemma}\label{ber}
Writing $\c^{\spli}=C_L$, we have
$$
\Ber(\c)^{\red}=L
$$
\end{lemma}
\begin{proof}
We know that $\Ber(\c)$ is isomorphic to $\D^{-1}$. In local super conformal co-ordinates, $\D$ is generated by $\frac{\partial}{\partial \theta}+\theta\frac{\partial}{\partial z}$. Recall that $\theta$ is an odd co-ordinate, so a section of $L$, and $\frac{\partial}{\partial \theta}$ is a section of $L^{\vee}$. Setting $\theta=0$ we obtain $\D^{\red}=L^{\vee}$, so the conclusion. More details about super-conformal co-ordinates can be found in \cite{Witten2} and \cite[Lemma 3.1]{WD}.

\end{proof}

\begin{theorem}[Pluri-Canonical Model]\label{CanMod}
Let $$f\colon \c\to B$$ be a susy curve over a super-manifold $B$. Fix an integer $\nu\geq 3$. Then $f_*(\Ber(\c)^{\otimes \nu})$ is a locally free sheaf of $\O_B$-modules. If $\nu$ is even, the rank of $f_*\Ber(\c)$ is
$$(\nu-1)g-\nu+1 \, | \, (2\nu-1)g-2\nu+1 \; .$$
If $\nu$ is odd, the rank of $f_*\Ber(\c)$ is
$$(2\nu-1)g-2\nu+1\, |\,  (\nu-1)g-\nu+1 \; .$$

If the genus $g$ of $\c$ is at least $4$, then $\Ber(\c)^{\otimes \nu}$ is relatively very ample. If $g=3$ or $g=2$, we need respectively $\nu\geq4$ or $\nu \geq 5$ for the very ampleness.
\end{theorem}

\begin{proof}

In this proof we follow the notations of Section \ref{preliminaries}. To show that $f_*\Ber(\c)^{\nu}$ is locally free we apply Criterion \ref{loc_free} and Lemma \ref{ber}. We need to show that 
$$\begin{sis}
& f^{\red}_*L^{\nu} \textrm{ and } f^{\red}_*L^{\nu+1}  \; \textrm{are vector bundles on} \; B^{\red};\\
& R^1f^{\red}_*L^{\nu}=R^1f^{\red}_*L^{\nu+1}=0.\\
\end{sis}$$
The degree of $L$ being $g-1$, these conditions are true if $\nu\geq 3$, but fail for $\nu\leq 2$. The statement about the rank also follows from Criterion \ref{loc_free}.
\medskip

Let us show that the bundle if $f$-very ample. We apply Criterion \ref{very_ample}. To fix the notation, assume that $\nu$ is even. First, we remark that
$$R^1f_*[(\Ber(\c)^{\nu})^{\red}\otimes I_x\otimes I_y]=R^1f_*[L^{\nu}\otimes I_x\otimes I_y]$$
By the classical Serre duality theorem, the right hand side vanishes for $\nu\geq 3$ if $g\geq 4$, for $\nu\geq 4$ if $g=3$ and for $\nu\geq 5$ if $g=2$. The second condition of Criterion \ref{loc_free} equally follows, because $(\Ber(\c)^{\nu})^{\red}\otimes L=L^{\nu+1}$. The same computation applies when $\nu$ is odd.
\end{proof}

\begin{remark}\label{gen}
Theorem \ref{CanMod} holds in a more general set-up: we have never used the susy structure. What we need is just that the line bundle $L=\Ber(\c)^{\red}$ is of degree $g-1$. 
\end{remark}

\begin{remark}[Canonical model]\label{canonical} The line bundle $\Ber(\c)$ is not $f$-very ample for a generic curve. For example, if $B$ is a point, $X=C_L$ and $h^0(C,L)=0$, it is well known that the rank of $f_*\Ber(\c)=H^0(X,\Ber(\c))$ is $0|g$, see for instance \cite[Section 8]{Witten2}; in this case, $\Ber(\c)$ can not be very ample for dimensional reasons. This means that we can not speak about the canonical model of a susy curve, but just about its pluri-canonical models. 

\end{remark}

\end{subsection}

\end{section}
\begin{section}{Moduli space and Hilbert Scheme}\label{hilb}

This section is purely speculative, because a theory of stacks and Hilbert schemes in the super setting has not be fully developed yet.

\smallskip

The moduli space $\sM_g$ of genus $g$ Super Riemann surfaces can be defined with the usual machinery of stacks: it is the pseudo-functor mapping a smooth super variety $B$ to the groupoid of susy curves over $B$. In other words, for any super-manifold $B$, the $B$-points of $\sM_g$, usually denoted by $\sM_g(B)$, are the susy curves over $B$. See \cite{SuperStacks} for an introduction to stacks in the super setting.

\smallskip

To simplify the notations, take $g\geq 4$ and $\nu=3$. In Theorem \ref{CanMod}, we showed that $\Ber(\c)^{\otimes 3}$ embeds every susy curve in a super projective space of dimension $2g-1|5g-5$. We thus want to consider the Hilbert scheme $H_g$ parametrising smooth connected sub-varieties of $\P:=\P^{2g-1|5g-5}$ with appropriate degree, dimension $1|1$, and embedded with $\Ber(\c)^{\otimes 3}$.

The Hilbert scheme $H_g$ can be defined as a functor; we do not know if it can be represented by a projective super-scheme, as in the classical case. 

Theorem \ref{CanMod} means that we have a morphism
$$
\iota \colon \sM_g \to H_g/\Aut(\P)
$$
The group $\Aut(\P)$ is studied in \cite{FKAut} and \cite{Penkov}. 

\begin{theorem}\footnote{This result was well known to the experts, we thank R. Donagi for discussing this theorem with us}
The morphism $\iota$ defined above is a closed embedding.
\end{theorem}
\begin{proof}

Following \cite[Tag 04XV]{stacks-project}, we have to show that $\iota$ is a monomorphism and universally closed

To show that $\iota$ is a monomorphism, we have to show that for any super-manifold $B$, the morphism $\iota \colon \sM_g(B) \to H_g/\Aut(\P)(B)$ is injective. To this end, we have to show that, given a $1|1$ dimensional super manifold over $B$, the susy structure $\D$, if it exists, it is unique. This is shown in \cite[Proposition 1.7]{DeligneLetter}. Let us sketch a proof that is slightly closer to \cite{Witten2}.

First, we prove the result in the split case $C_L$. A split $1|1$ manifold is equivalent to a curve $C$ with a line bundle $L$; to have a susy structure $\D$ on $C$, $L$ must be a theta characteristic. The distribution $\D$ is uniquely determined by $L$ and vice versa, so the claim. 

Once the split case is established, we have to show that we can not deform $\D$ keeping the underlining super-manifold fixed. It is enough to prove it for a first order deformation. To do this we argue as follows. Let $\S$ be the sheaf of super-conformal vector field; that is the vector field that commute with $\D$. A deformation of a susy curve is given by an element of $H^1(C,\S)$. We can forget the susy structure and see the deformation just as a deformation of a $1|1$ super-mainifold. This corresponds to look at the map
$$
i \colon H^1(C,\S)\to H^1(C,TC)
$$
induced by the inclusion of $\S$ in $TC$. The claim is equivalent to the injectivity of $i$; indeed, this means that if we deform the susy structure then we deform also the underlying super-manifold. As explained in \cite[Page 11]{Witten2}, the projection
$$
\pi \colon TC\to TC/\D\cong \D^2
$$
induces an isomorphism between $\S$ and $\D^2$. In other words, in the (non-exact) sequence
$$
\S \xrightarrow{i} TC \xrightarrow{\pi} \D^2\cong\S
$$
we have $i \circ \pi= Id$. If we take to cohomology, we have
$$H^1(C,\S) \xrightarrow{i} H^1(C,TC) \xrightarrow{\pi} H^1(C,\S)$$
with again $i \circ \pi= Id$, so $i$ is injective at the level of $H^1$.

\medskip

To prove that $\iota$ is universally closed, we show that the image is the fixed locus of an involution of the co-domain; to this end, we introduce duality. Given a $1|1$ super-manifold $X$ we can define the dual curve $X^{\vee}$. There are a few ways to define $X^{\vee}$. A possibility is the moduli space of $0|1$ dimensional sub-manifold of $X$ whose intersection with the reduced space is a point. For more details see e.g. \cite[Section 9]{Witten2} or \cite[Sections 2.2 and 2.3]{Ber}.

Because of Remark \ref{gen}, we can see  $H_g/\Aut(\P)$ as the moduli space of dimension $1|1$ super-manifolds  such that $\deg \Ber(X)^{\red}=g-1$. The duality preserves these invariants, in particular $\deg \Ber(X^{\vee})^{\red}=2g-2-\deg \Ber(X)^{\red}=g-1$, cf. \cite[Example 2.2.3]{Ber}. We conclude that the duality gives an order $2$ automorphism of $H_g$. Susy curves are exactly the curves which are auto-dual, cf. \cite[Section 2.3]{Ber} or \cite[Section 9]{Witten2}; this means that $\sM_g$ is the fixed locus of this automorphism of $H_g/\Aut(\P)$.

\end{proof}

\end{section}

\end{document}